\documentclass[12pt,a4]{article}
\usepackage{color}
\usepackage{url}
\usepackage{latexsym}
\usepackage{amsmath,amssymb}
\usepackage[round,sort]{natbib}

\begin{document}
\title{\large \bf{On the Transition Period of Implementing New Mathematics Curriculum for Foundation Engineering Students}}%
\author{\normalsize N. Karjanto$^{\tiny 1}$ and R. Osman$^{\small 2}$\\%
{\footnotesize $^{1}$Department Mathematics, School of Science and Technology}\\
{\footnotesize Nazarbayev University, Astana, KAZAKHSTAN}\\
{\footnotesize{\url{natanael.karjanto@nu.edu.kz}}}\\
{\footnotesize and}\\
{\footnotesize Department of Mathematics, University College}\\
{\footnotesize Sungkyunkwan University, Natural Science Campus}\\
{\footnotesize 300 Cheoncheon-dong, Jangan-gu, Suwon, SOUTH KOREA}\\
{\footnotesize{\url{natanael@skku.edu}}}\\
{\footnotesize $^{2}$Department of Applied Mathematics, Faculty of Engineering}\\
{\footnotesize The University of Nottingham Malaysia Campus}\\
{\footnotesize{Jalan Broga, Semenyih 43500, Selangor Darul Ehsan, MALAYSIA}}\\%
}%
\date{}
\maketitle\thispagestyle{empty}

\begin{abstract}
\noindent \footnotesize{An overview on several mathematics modules in the transition period of introducing a  new curriculum for the Foundation programme in Engineering at the University of Nottingham Malaysia Campus is discussed in this paper. In order to progress to Undergraduate programmes in Engineering, previously the students must complete three mathematics modules of 40 credit points in total, for which one of them {\color{black} was} a year-long module with 20 credit points. Currently under the new curriculum, the students are required to complete five mathematics modules with 10 credit points each. The new curriculum gives positive impacts for both the lecturers and the students in {\color{black} terms} of material organization, fully utilizing textbooks and a new arrangement for tutorial sessions. The new curriculum also provides the students with stronger mathematical background in critical thinking and problem solving skills to equip them to embark the Undergraduate programmes in  Engineering.}\bigskip

\end{abstract}

\section{Introduction}

\subsection{\color{black} Background and motivation}
The University of Nottingham in Malaysia is one of the branch campusses of the University of Nottingham in the UK. Another overseas branch campus of the university is located in Ningbo, a seaport city of northeastern Zhejiang Province, eastern part of People's Republic of China. Our campus houses four {\color{black} Faculties} and the Department of Applied Mathematics belongs to the Faculty of Engineering. Although our department does not really offer any study programme majoring in Mathematics, it plays an important role as a service department not only for {\color{black} the} Undergraduate programme but also for {\color{black} the} Foundation programme in Engineering. Currently, there are four engineering departments within the Faculty of Engineering that are being served by our department for their mathematics modules: Electrical and Electronic Engineering, Chemical and Environmental Engineering, Mechanical, Materials and Manufacturing Engineering and Civil Engineering.

In contrast to the North American system where the students will normally spend four years in their undergraduate study, the British system of higher education is adopted in our curricula, where the undergraduate level is divided into three years~\citep{framework}. This means that the students who enter to Year~1 have an entry requirement of A~(advanced) level or a similar qualification. However, a significant number of our potential students possess only an O~(ordinary) level qualification, particularly the local Malaysian students who have completed their SPM ({\sl Sijil Pelajaran Malaysia} or the Malaysian Certificate of Education). This SPM is equivalent to the British GCSE (General Certificate of Secondary Education), and provides the opportunity for the students to continue their {\color{black} studies} to pre-university (Foundation) level. For those who are interested to obtain a better grasp on the education system in Malaysia, including the status of mathematics teaching and learning in the country, the reader is strongly encouraged to consult an article by~\cite{Sam09}.

In order to accommodate the students with an O~level qualification, the {\color{black} Faculty} made a special arrangement, known as the Foundation programme in Engineering. Similar Foundation programmes in Computer Science, Bioscience and Business are also offered by the other {\color{black} Faculties} in the branch campus. The Foundation programme in Engineering lasts for three semesters which normally begins either in April or in July and {\color{black} is completed} in May of the following year. The aim of this programme is to {\color{black} support} the students with merely an O~level qualification to enter the Undergraduate programme in Engineering. In this case, at the end of the third semester of their Foundation study, the students would have a similar level to those with an A~level qualification. After completing the Foundation programme and satisfying the progression rules, the students may proceed to choose one study programme offered by the four engineering departments. The Foundation programme in Engineering covers a set of modules that contains {\color{black} topics} on Mathematics, Physical Sciences (including Physics and Chemistry), Information Technology, English and some other modules for the local needs.

{\color{black} In this paper, the curricula for mathematics modules at the Foundation level will be discussed. The purpose of this article   is not only to provide an overview of these mathematics curricula but also to explain some pedagogical and content {\color{black} thoughts} behind the need to adjust and to modify the curricula. The paper {\color{black} also} presents some problems associated with {\color{black} the} teaching {\color{black} of} a number {\color{black} of} mathematics modules to students enrolled in the Foundation programme.
Explicit connections between our new curriculum and current development in mathematics curriculum are discussed in this paper. In a more general context, we show that by implementing the new curriculum, some major problems faced in teaching mathematics at the early undergraduate (Foundation) level for Engineering programme could be solved. 
The solutions proposed are in line with current {\color{black} thoughts on} teaching and learning mathematics in general. Although empirical results might not be readily available, some anecdotal evidences and initial outcomes {\color{black} from} the implementation of the new curriculum are also included.}

Since the focus of this paper is specifically on mathematics curriculum for mathematics modules at the Foundation level, a spacious room is absolutely still available for discussion of mathematics curriculum at the undergraduate level for Engineering programme in general or even for other study programmes that need mathematics.

\subsection{\color{black} Literature review}

A general overview on problem solving in the mathematics curriculum is reported by~\cite{Schoenfeld83}.
{\color{black}
The four phases of `Multidimensional Problem-Solving Framework' (MPSF), i.e. orientation, planning, executing and checking, that characterizes various problem solving attributes are described by~\cite{Carlson05}. Other authors proposed five-phased model of engagement, transformation-formulation, implementation, evaluation and internalization to describe (meta)cognitive approaches for problem solving~\citep{Yimer10}. A theory of `goal-oriented decision-making' in mathematics problem solving is proposed by~\cite{Schoenfeld10}. Recently, an approach of `grounded theory' for problem solving skill in engineering problems has been published by~\cite{Harlim13}.
}

A research on students' understanding of functions and its importance for the undergraduate mathematics curriculum is discussed by~\cite{Thompson94}. The design of the mathematics curriculum for engineers in an Australian university has been addressed by~\cite{Varsavsky95}. Her findings indicate that to achieve cohesion and to make the course more meaningful to students, the design of the mathematics curriculum must be done in close collaboration between the mathematics department and the engineering faculty. A reform in undergraduate mathematics curriculum with more emphasis on social and pedagogical skills is presented by~\cite{Pesonen00}. The readers who are interested in the history of mathematics curriculum, the culture of mathematics, gender, and social justice issues in mathematics may consult an article written by~\cite{Teese00}.

An examination of the various forces which act on mathematics curriculum and on curriculum trends in both at local and US national level is presented by~\cite{Hillel02}. The article covers some issues including undergraduate programmes, specific courses, mathematical content, degree of rigour, modes of delivery and interaction and assessment schemes. An article discussing the importance of mathematics in the university education for engineers is reported by~\cite{Kent03}. The authors emphasize mathematical skills, studio-based and problem based learning techniques and use of information technology (IT) in teaching mathematics to engineering students. An identification and enhancement of mathematical understanding among engineering undergraduate students at MIT by improving mathematics curriculum is presented by~\cite{Willcox04}. A report describing research to create explicit links between engineering courses and upstream mathematics courses at MIT is discussed by~\cite{Allaire04}.

{\color{black}
A suggestion for a theoretical model based on the anthropological notion of a modern-day rite of passage for a transition from secondary school to university mathematics in the context of North-American educational context is featured by~\cite{Clark08}. An enriched version of the model, with the added notions of cognitive conflict and culture shock, has been addressed by the same authors~\citep{Clark09}. Their model suggests that the transition from high school to university mathematics requires a proper environment and there is a need for filling in the temporal gap between the two stages with a set of meaningful activities.

The transition of mathematics learning from secondary school to university from the students' perspective has been addressed amongst others by~\cite{Barnard03} and \cite{Hernandez11}.
A large body of literature exists discussing surrounding issues of the transition of teaching and learning in mathematics from secondary to tertiary levels, amongst others are~\citep{Wood01,Jourdan07,Brandell08,Brown08}.
Some findings from a project analyzing the transition from secondary to tertiary education in mathematics from the teachers' and lecturers' perspectives is reported by~\cite{Hong09}. The results provide evidence that each group lacks a clear understanding of the issues involved in the transition from the other's perspective and there is an urgent need to improve the communication between the two parties.}
An excellent article on the changes in thinking involved in the transition from school mathematics to formal proof in pure mathematics at university {\color{black} has been written} by~\cite{Tall08}. It is interesting {\color{black} for} the students enrolled in the Foundation level {\color{black} to} experience this transition level from high school mathematics to university mathematics.


This paper {\color{black} fills the literature gap of the transition level from secondary mathematics education to undergraduate mathematics education and} is organized as follows. After this introduction, Section~\ref{issue} addresses some problems and issues in teaching mathematics modules at the Foundation levels. Additionally, Section~\ref{current} discusses current thinking on mathematics teaching and learning, in particular the curriculum which emphasize the students acquiring problem solving skills.  Furthermore, Section~\ref{curricula} gives an overview on the old and the new curricula of mathematics modules at the Foundation level, including some changes that have been implemented. Next, Section~\ref{reaction} provides some initial reactions to the implementation of this new curriculum. Specific descriptions addressing the concerns are also discussed in this section. Finally, Section~\ref{conclusion} gives {\color{black} the} conclusion {\color{black} of} our discussion.

\section{Issues surrounding teaching Foundation mathematics} \label{issue}

It is observed that many students who are currently enrolled in the Undergraduate programmes in Engineering possess weak background in basic mathematics. This observation is based on our experience in teaching several mathematics modules to the Undergraduate students. Many of the students have forgotten some basic mathematical concepts that were introduced to them earlier when they were enrolled in the Foundation programme. For instance, some are confused to conclude that when $f(x - 2) = 1/x$, then $f(x) = 1/(x + 2)$. Many find also difficult to conclude that if ${\displaystyle g(t - 2) = \frac{t}{(t - 2)^{2} + 1}}$, then ${\displaystyle g(t) = \frac{t + 2}{t^{2} + 1}}$. It can even be surprising that a Year 3 student was not able to perform simple integration and derivation, such as evaluating ${\displaystyle \frac{{\rm d}}{{\rm d}x} \int_{1}^{x} 2s \, {\rm d}s}$. Actually this is simply the First Fundamental Theorem of Calculus, which is taught at the Foundation level{\color{black}; see~\cite{Smith08} or any other books on Calculus.}

In addition, many students lack of ability in critical thinking and problem solving skills. Even though the enrollment to the programme has been through a selection process, we may still enrolled the students who are rather weak in mathematics but quite strong in other fields and thus the current situation is simply inevitable. The cr\`{e}me de la cr\`{e}me students in the country prefer to choose {\color{black} well-known local universities} instead of our private institution, mainly for financial {\color{black} reasons}. A similar situation occurs for our international students. Even though we manage to capture excellent students, the indisputable fact is that on average, many of our students are academically at the intermediate level, particularly in mathematical ability. It is our responsibility to train them to think critically and to possess problem solving {\color{black} skills}.

Furthermore, it is also observed that a year-long module of 20 credit points is always tougher for the students than two separate modules of 10 credit points each. In the context of our university curriculum, a~10 credit {\color{black} points'} module would provide the students a two-hour lecture and a one-hour tutorial sessions each week.
A final exam covering the semester's material is conducted at the end of each semester. For a year-long module with 20 credit points, a similar class session with the one from the 10 credit points module is also adopted. However, there is no final exam at the end of the first semester; the final exam is conducted at the end of the second semester {\color{black} and} covers the material from both semesters. This particular type of arrangement {\color{black} imposes} a heavy burden to the students since they have to cover two-semester materials for the final exam.

With this in mind, we have proposed to split the 20 credit points modules into two modules with 10 credit points each. As a consequence, there will be two final exams at the end of the first and the second semester, respectively. Under this new arrangement, the burden for the students will be alleviated since they {\color{black} simply need} to prepare the material for each semester only. Moreover, in response to improve students' problem solving {\color{black} skills} and their critical thinking ability, we have also developed our curriculum {\color{black} with respect to these} needs. This strategy is in line with the current {\color{black} thoughts on} teaching and learning mathematics in general and engineering mathematics in particular, as {\color{black} it will be} discussed in the following section.

\section{Current {\color{black} thoughts on} mathematics teaching and\\ learning} \label{current}

The main aim of implementing this new curriculum is to improve the confidence, mathematical knowledge and fluency of the students in problem solving. The expected learning outcome consists of four important components: knowledge and understanding, intellectual skills, professional skills and transferable skills. Basically, the first component will depend {\color{black} on the} contents of the module. However, the latter three components can be described in {\color{black} a} more general framework as follows.
\begin{itemize}
\item {\color{black} Concerning the} intellectual skills, the students who complete the module should be able to reason logically and work analytically, perform high levels of accuracy, manipulate mathematical formulas, algebraic equations and standard functions and apply fundamental mathematical concepts to {\color{black} routine problems} in engineering or science.

\item {\color{black} Concerning the} professional skills, the students who complete the module should be able to construct and present mathematical arguments with accuracy and clarity as well as apply basic solution techniques learned to mathematical problems arising in the study of engineering or science.

\item Finally, the students who complete the module are also expected to obtain transferable skills, i.e. to communicate mathematical arguments using standard terminology, express {\color{black} the} ideas of solution of mathematical problems appropriately and effectively and use an integrated software package to enhance learning and practice their problem solving skills.
\end{itemize}

It is important to note that our new curriculum heavily emphasizes on teaching and learning through problem solving.
Problem sheets are distributed each week and a number of selected questions are discussed in the following week during the tutorial session. The type of {\color{black} questions} in the problem sheets ranges from simple to more challenging {\color{black} ones}.

The implemented teaching and tutorial sessions still rely heavily on {\color{black} the} traditional method, where the instructor tends to dominate the entire teaching session. Although the students are also strongly encouraged to participate actively during the tutorial session, it remains a challenge to have an interactive teaching session, since many students are timid by nature.

The third component of {\color{black} our} theoretical framework is the {\color{black} process} of gathering and analyzing data. For this, we still yet have to wait for the {\color{black} outcomes} of implementing the new curriculum. However, the previous data of the students' {\color{black} results}, as well as {\color{black} the} class observation of the students' {\color{black} abilities, suggest that} there is an urgent need to improve the curriculum and teaching emphasis.

It is interesting to note that our new curriculum is developed in line with current {\color{black} trends} in teaching mathematics in general, as we {\color{black} have} discovered in the literature. A description of a particular framework for research and curriculum development in undergraduate mathematics education with some examples of its application is given by~\cite{Asiala96}. The authors describe certain mental {\color{black} structures} for learning mathematics, including actions, processes, objects, and schemas, {\color{black} also known as the acronym APOS,} and the relationships among these constructions. The components of the ACE teaching cycle (activities, class discussion, and exercises), cooperative learning and the use of a mathematical programming language are also explained thoroughly.

{\color{black}
The theory of APOS based on Piaget's principle that an individual learn mathematics by applying particular mental mechanism to build specific mental structures and uses these to solve problems related to particular situations~\citep{Piaget70}. The main mechanisms are interiorization and encapsulation and the related structures are APOS themselves. There are three stages for the implementation of APOS theory as a framework for teaching and learning mathematics, i.e. theoretical analysis, instructional sequences and data collection and analysis~\citep{Dubinsky01}. Based on this APOS theory, the pedagogical approach of ACE teaching cycle that encourages active student learning is developed. Although the activities meant by the previous authors are computer related activities~\citep{Asiala96,Voskoglou13}, in our new curriculum, the activities also include the traditional classroom activities, where the instructor poses some mathematical problems to the students before explaining and covering a particular topic~\citep{Maharaj10}.

The following provides an example mentioned earlier in the context of APOS theory. {\bf Action}: A student who wants to find an expression of the function $f(x)$ given the form $f(x - a) = 1/x$, $a \neq 0$, $x \in \mathbb{R}$ or finding $g(t)$ given a rational function $g(t - b) = h_1(t)/h_2(t - b)$, where $h_1$ and $h_2$ are other functions in $t$ and $b \neq 0$, $t \in \mathbb{R}$, can find the expression for the functions $f$ and $g$ directly by observing the expression at the right-hand sides. {\bf Process}: an individual with a process of understanding the substitution will construct a mental process by replacing old variables $x$ and $t$ with new variables $\xi = x - a$ and $\tau = t - b$ and rewrite the new variables $\xi$ and $\tau$ with the old ones $x$ and $t$ once the process is complete. {\bf Object}: the student could think those functions as geometrical object and changing the variables mean shifting the graph horizontally to the left $(a, \ b > 0)$ or to the right $(a, \ b < 0)$. {\bf Schema}: the individual organizes all the other three components of actions, processes and objects into a coherent framework, thus a complete understanding of problem solving is attained.
}

The three components of the theoretical framework for mathematics education mentioned in~\citep{Asiala96} are adopted implicitly in our new curriculum. The first component, i.e. theoretical analysis, where the students acquire knowledge and understanding, is covered in the first learning outcome of the new curriculum. The second component of the theoretical framework is implementation of instruction. This postulates certain specific mental construction that the instruction should foster. In connection to our new curriculum, the instruction of the new curriculum is tailored to the desired learning outcome, which includes the intellectual, professional and transferable skills. The implemented instruction method will help the students to use certain constructions in different situation and to develop problem solving skills. Finally, the third component of the theoretical framework is the collection and analysis of data. Although at this stage the outcome of the new curriculum is not fully documented yet, some anecdotal evidences and initial reactions, however, are discussed in Section~\ref{reaction}.

Another {\color{black} work} emphasizes problem solving as a basis for reform in mathematics curriculum and instruction~\citep{Hiebert96}. The authors discuss the history of problem solving in the curriculum that has been infused with a distinction between acquiring knowledge and applying it. They propose an alternative principle by building on the idea of reflective inquiry, arguing that the approach would facilitate students' understanding. A number of mathematics curricula for engineers have also been designed to tailor the need for critical thinking and problem solving skills. See for instance amongst others~\citep{Lichtenberger02,Lesh03,Gainsburg06,Hurford09} and the references therein.

The problem solving approach proposed by~\cite{Hiebert96} allows students to wonder why things are {\color{black} so}, to search for solutions and to resolve incongruities. Our new curriculum is also designed in a similar line of ideas. The instruction is tailored so that the students are engaged in problem solving activities and thus {\color{black} being trained to} reason logically and to think critically (intellectual skills). Historically, the mathematics curriculum in general has been shaped by concerns about preparation for the workplace and for life outside of school~\citep{Stanic88}. Problem solving has been used as a vehicle to reach this goal. The professional skills of our new curriculum include some aspects of students' construction when they encounter problems in science and engineering, later during their study of even after they leave university. Finally, the transferable skills {\color{black} emphasized by}~\cite{Hiebert96}, {\color{black} who propose} an implementation of problem-based learning or case-based instruction, {\color{black} are also adopted in the curriculum}.

{\color{black}
More recent views on problem solving has been mentioned in the introduction of this paper. New insights on problem-solving process by offering multidimensional framework to investigate, analyze and explain mathematical behaviour are described by~\cite{Carlson05}. The authors explains four phases of `Multidimensional Problem-Solving Framework' (MPSF) of orientation, planning, executing and checking, where various problem-solving attributes, including their roles and significance during each phases, are characterized by MPSF.

An excellent scholarly work on the understanding of how and why pedagogical decision-making happens in the course of teaching and learning can be found in~\citep{Schoenfeld10}. The author proposes a theory of `goal-oriented decision-making' and provides particular examples of problem solving in mathematics where teachers' instructional decision making is a function of instruction goals, resources and orientations. A five-phased model to describe the range of cognitive and metacognitive approaches used for mathematical problem solving is proposed by~\cite{Yimer10}. The five-phases are engagement, transformation-formulation, implementation, evaluation and internalization. Findings on how engineers develop their problem solving skills with implications on general educational strategy, the development and the implementation of computer technology for engineering problem solving are published recently~\citep{Harlim13}. The authors utilized an exploratory approach of `grounded theory' to understand the complexities of engineering problem solving.
}



The following section explains the differences between the old and the new curricula and a new tutorial arrangement conducted after implementing the new curriculum.

\section{The old and new mathematics curricula} \label{curricula}

\subsection{The old curriculum}

The {\color{black} old} curriculum referred in this paper is the one that has been implemented until the academic year 2008/2009. Under this former curriculum, three mathematics modules with total of credit points of 40 are given to the Foundation students. These are HG1BMT Basic Mathematical Techniques (10 credit points) offered in Semester~0, HG1FND Foundation Mathematics (20 credit points) offered in Semesters~1~and~2 and HG1M02 Applied Algebra (10 credit points) offered in Semester~2. Thus, the students who have completed their Foundation programme have taken three mathematics modules with total of 40 credit points.

The {\color{black} following} are the {\color{black} summaries} of the content of each {\color{black} of the} mathematics modules under the former curriculum.
\begin{itemize}
  \item HG1BMT Basic Mathematical Techniques\\
  This module provides a basic course in algebra and {\color{black} introduces} some basic knowledge on functions and analytic geometry. This module {\color{black} cover basic} algebra, inequalities, polynomials, functions and graphs, coordinate geometry, conic sections, sequences and series, binomial expansion and {\color{black} partial fraction decomposition}.

  \item HG1FND Foundation Mathematics\\
  This module provides the basic topics {\color{black} of} differential and integral calculus. It covers trigonometric functions, complex numbers, differentiation, {\color{black} applications of derivatives, integrals}, numerical integration, curve sketching and the binomial theorem for any rational index.

  \item HG1M02 Applied Algebra\\
  In this module, some of the fundamental concepts of vector and linear algebra that arise naturally in many engineering circumstances are introduced. Problems that demonstrate the applicability of the theory are covered and this enables students to develop a facility at applying the techniques for themselves. This is considered to be a crucial aspect of the training of a modern engineer. This module {\color{black} covers} determinants, vector algebra and its applications to three-dimensional problems in geometry, vector differential operators, matrices and systems of equations.
\end{itemize}

\subsection{The new curriculum}

The new mathematics curriculum {\color{black} is} the one that is {\color{black} currently} being implemented from the academic year 2009/2010 onward, starting since April 2009. Under this new curriculum, there is no significant change to HG1M02 Applied Algebra. However, the 20 credit {\color{black} points'} module HG1FND Foundation Mathematics is split into two 10 credit points modules, to become F40CA1 Calculus~1 and F40CA2 Calculus~2, with some slight variations. These modules are offered in Semester~1 and Semester~2, respectively and basically contain differential and integral calculus. A 10 credit {\color{black} points'} module which is similar to HG1BMT Basic Mathematical Techniques is brought in under a new name: F40FNA Foundation Algebra, which is offered in Semester~0. Furthermore, a new 10 credit points module is introduced under the new curriculum and is offered in Semester~1, namely F40FMT Mathematical Techniques.

The {\color{black} following are the summaries of the contents of each mathematics module} under the new curriculum.
\begin{itemize}
  \item F40FNA Foundation Algebra\\
  This module is offered in Semester 0 and provides a basic course in algebra and trigonometry. It also introduces the students to skills in core mathematical techniques. This module will cover indices and rules of algebra, quadratic equations, polynomials, {\color{black} partial fraction decomposition}, trigonometry, {\color{black} logarithms}, inequalities, sequences, series and binomial expansion.

  \item F40FMT Mathematical Techniques\\
  This module is offered in Semester 1. It introduces the complex number system and it provides a basic course in elementary statistics.
  Algebraic manipulation and operations on complex numbers are introduced with the aid of an Argand diagram in the complex plane. Initial key elements of definition, manipulation and graphical representation of data are introduced prior to establishing statistical techniques used in the analysis of problems in engineering and physical sciences. Application to solving real life problems is developed. The module will cover complex numbers, basic set theory, graphical representation of data, numerical descriptive, probability and counting techniques, discrete probability distribution (binomial distribution) and continuous probability distribution (normal distribution).

  \item F40CA1 Calculus~1\\
  This module is also offered in Semester 1 and provides a basic course in differential and integral calculus. Initial key elements of definition, manipulation and graphical representation of functions are introduced prior to establishing techniques of calculus used in the analysis of problems in engineering and physical sciences. {\color{black} Applications in} solving real life problems {\color{black} are} developed. The module {\color{black} covers} functions and graphs, limits and continuity, techniques of differentiation, applications of differentiation and curve sketching.

  \item F40CA2 Calculus~2\\
  This module is offered in Semester 2 and provides a basic course in integral calculus. Initially, formulas for integration as an antiderivative of functions are introduced prior to establishing techniques of integration for more complicated functions. For {\color{black} applications}, integrals are used to evaluate the area under a graph, the volume of revolution, mean value, and root mean square value of a function. For a non integrable function, numerical integration is introduced. The technique of the Newton-Raphson method is also demonstrated. Finally, the conic {\color{black} sections} with their equations and the {\color{black} parametric equations} of curves are explained. This module {\color{black} covers} indefinite and definite integrals by formulas,  techniques of integration, integration by parts, {\color{black} applications of the} integral, differential equations with separable variables, numerical integration, conic {\color{black} sections} and parametric equation of curves.

  \item HG1M02 Applied Algebra\\
  This module is also offered in Semester 2 and provides the basic elements of vector algebra and linear algebra and their applications to simple engineering situations. It  introduces the modeling of basic engineering situations in terms of multi-dimensional models. Initially the key elements of definitions and manipulations of basic mathematical skills and mathematical techniques in matrices and vectors are introduced prior to modeling and analyzing problems related to engineering situations. The module {\color{black} covers} vectors, matrices, system of equations, Gauss and Gauss-Jordan elimination and Cramer's rule.
\end{itemize}

We observe that trigonometry and {\color{black} logarithms} which formerly {\color{black} were presented} in Semester~1, they {\color{black} are presented} in Semester~0 under the new curriculum. Curve sketching moves from Semester~2 into Semester~1. Complex number stays in Semester~1, but it is taught under the new module Mathematical Techniques. Furthermore, this new module simply contains new materials that are not really {\color{black} covered} in the old curriculum. These are basic set theory and introduction to probability and statistics. Introduction to differential equations with separable {\color{black} variables is a new topic} and is {\color{black} covered} in Semester~2 under Calculus~2. Basically, there is no significant change between the former and the new curricula for HG1M02 Applied Algebra.

\subsection{A new tutorial arrangement}

A new arrangement on tutorial session has been introduced starting in the current academic year for the Foundation students enrolled in April 2009 intake. This arrangement is particularly beneficial for those who are less strong in the mathematical competency compared to the average of the class within the same batch. A scholastic test is given at the beginning of the semester to distinguish the weaker group of students from the stronger one. Thus, instead of following the $2 + 1$ sessions, the weaker students will follow $3 + 1$ sessions. It means that they will receive a three-hour lecture and a one-hour tutorial session. Since a measure of flexibility is allowed, the lecturer might implement $2 + 2$ sessions, where the students will have more opportunity to do more problem solving and exercises during the class, since a two-hour session is devoted for the tutorial. After implementing this arrangement, it is expected that at the end of the semester, the gap of the mathematical competency between the two groups can be decreased. We are eagerly awaiting and curiously anticipating to the success of this novel idea.

\section{{\color{black} Research methodology} and its initial findings} \label{reaction}

This section explains the research methodology implemented in {\color{black}} study. It also reports some anecdotal evidences and initial findings {\color{black} concerning} the implementation of new curriculum from both the instructors' and students' viewpoints.

\subsection{\color{black} Research methodology}
{\color{black} A methodology of simple qualitative research interview} is implemented in this study, {\color{black} without explicit questionnaire and realizes upon discussions and arguments}. {\color{black} Since there is a wide variation among interviewing approaches, the approach of `unstructured interview' is adopted.
Although no interview can truly be considered unstructured, what is meant by unstructured interview in this context is more or less equivalent to a guided conversation, where the interview is conducted in conjunction with the collection of observation data. This is different with semi-structured interview where the interview is conducted as the sole data source for a qualitative research project and usually clinical type in nature.
Information on different formats of qualitative interviews and qualitative methods used in mathematics education research has been explored by~\cite{Romberg92}, \cite{Schoenfeld94} and \cite{DiCicco06}.}

We interviewed six faculty members of the Department of Applied Mathematics at UNMC. All of these colleagues have experience in teaching a number of mathematics modules at the Foundation levels using the old curriculum. Out of six faculty members, four of them have experience teaching mathematics modules using the new curriculum. Two of the faculty members are not interviewed since they have left the faculty to continue their PhDs but they have contributed to the development of the new curriculum.

We also interviewed the students who have completed their Foundation programme and the students who are still enrolled in the Foundation programme. The former cohort of students experiences the old mathematics curriculum but does not experience the new curriculum. Conversely, the latter cohort of students experiences the new curriculum but does not have any idea the situation for the old curriculum. For each cohort, we explain the differences in the curriculum and some changes that we have implemented. Since the students are generally a little bit intimidated when their instructors interviewed them, they {\color{black} were} interviewed in a rather informal settings, i.e. {\color{black} during office hour,} during lunch break, during class intermission and during the study week period before the final exam.

{\color{black} The following are several examples of the questions during the guided conversation of the unstructured interviews. What are some differences between the old and the new curricula? What aspects make the new curriculum more suitable for faculty and for students? Why is splitting a 20 credit points' module to two-10 credit points' modules beneficial? What are the benefits of a new arrangement for tutorial sessions? Why adopting a particular textbook for certain course is very practical? What do you think regarding the overall impression of the new curriculum?}
Some anecdotal evidences and initial reactions regarding the implementation of the new curriculum {\color{black} is reported in the following subsection}.

{\color{black} \subsection{Initial reactions}}

All members of academic staff {\color{black} of} the department  have responded positively to the content, topical arrangement and textbook choice implemented in the new curriculum. Formerly, the interdependence between some topics covered in the old curriculum is rather clumsy. For instance, topics on trigonometry and complex number are presented in the earlier part of HG1FND Foundation Mathematics and binomial expansion for any rational index is presented at the end of the module. In fact, the content core of Foundation Mathematics module is Calculus. Thus, although these three topics might be essential for students' understanding in mathematics, they may distract the main purpose for the module, i.e. to cover Calculus. Furthermore, having these three topics included has made it very challenging for us to {\color{black} select a suitable} textbook. Although many Calculus textbooks cover these topics in the appendix or as a quick review to remind the readers, they are not as thorough as we would like to be. For instance, a thorough discussion on trigonometry and trigonometric functions are generally found in the textbooks that also discuss Algebra, for instance {\sl Algebra and Trigonometry} by~\cite{Sullivan11}. The reason is that trigonometry is considered as a prerequisite for any general course in Calculus.

On the other hand, the topics covered from each module in the new curriculum are more specific and have strong interrelation. As instructors, we are very glad that some topics have been moved to another module. The coverage for Calculus has also been extended. For instance, trigonometry and binomial expansion have been transferred to F40FNA Foundation Algebra. Complex number is given in the new module F40FMT Mathematical Techniques. The division of the materials between F40CA1 Calculus 1 and F40CA2 Calculus 2 is quite distinct. Basically, Calculus 1 covers differentiation, {\color{black} its applications} and all related things to it and Calculus 2 covers integration and its accompanying techniques and applications. From {\color{black} the} instructors' perspective, we affirm that the new curriculum has a better organization {\color{black} with respect to the topics} distribution.

As instructors, we also respond positively to the admonition of implementing textbooks in our teaching. In the old curriculum, there is no particular textbook that neither we and the students could use due to rather peculiar topical arrangement. As a consequence, we spend more time in compiling lecture notes for the students, with inevitable minor typographical errors here and there. This definitely {\color{black} took} our time for teaching preparation and in turns may {\color{black} affected} the quality of our teaching. However, in the new curriculum, some particular textbooks have been chosen as main sources of reference for the corresponding modules. For example, we adopted an A-level mathematics textbook from UK for our Foundation Algebra module, written by~\cite{Bostock02}.

Adopting and implementing textbooks for mathematics modules will not only help us in our teaching preparation but also it saves us uncountable precious time in preparing lecture notes as well as preparing sets of problem sheet. {\color{black} Consequently, a high quality teaching is delivered every session. In addition,} both teaching and tutorial sessions are much easier to handle since we simply refer to the textbook instead of heavily depended on our self-built lecture notes. In addition, when assigning certain exercise questions to the students, we simply refer to the textbook instead of rewriting in sheet of papers and then distribute them to the students. Since textbooks have gone to editorial process, the number of errors is usually very minimal. This is another advantage of adopting and implementing textbooks as part of the new curriculum. Overall, the new curriculum has {\color{black} saved} us time both from teaching preparation and from administrative aspects.

The students have expressed their positive responses toward the implementation of the new curriculum. The new students who are enrolled in the Foundation programme are high school graduates and might not be aware regarding the change of the curriculum at the university. However, the students who have completed their Foundation programme and enrolled in the Undergraduate programme have discovered that the new curriculum is better {\color{black} than the old one}. They respond positively in term of re-arrangement of the topics, an added extra module and an extra hour for class interaction. They wish that that could enjoy the privilege and they affirmatively confirm that the new curriculum is more favorable for their younger peers. These comments particularly come from the older students who are rather weak in their mathematical and analytical skills, but have a little opportunity to improve them and yet they are willing to learn. Provided that these students put sufficient effort to study, the new curriculum contributes in helping them to improve their study skills.

Furthermore, even though the new students might not be aware regarding the difference between the old and the new curriculum, many still respond enthusiastically to the new curriculum. One particular positive response is to the arrangement of an extra hour class interaction. This extra hour might be use either for teaching or for tutorial session or both, depending on the need arises in the classroom. The instructor in this case has a measure of flexibility to implement which session is best suited for {\color{black} the need} his/her class. The students who are rather weak in mathematics would appreciate since the teaching is paced appropriately. The students are also beneficial when a step by step explanation is given, rather than jumping some steps. Importantly, the instructor could give guidance in problem solving and the students have more opportunity to ask and discuss the problems. Interestingly, some students who are qualified for a normal session have expressed their desire to join their peers on the slower session with an extra contact class hour. Although initially we are a little bit hesitant to arrange two groups with different number of hours for class interaction, the result is quite the opposite. We also observed that the students are working more diligently on the problem solving session and class discussion. So, the overall impression from the students' perspective regarding the new curriculum is also very positive and encouraging.\\

\section{Conclusion} \label{conclusion}

A brief overview on mathematics curriculum for students enrolled in Foundation programme in Engineering at the University of Nottingham Malaysia Campus has been presented in this paper. The relationship between the new curriculum to current {\color{black} trends} in teaching mathematics in general has also been discussed. It is observed that under the new curriculum, more materials are given to the students compared to the previous one. The purpose of this new curriculum is to {\color{black} increase the students'} mathematical ability {\color{black} in order to be} able to thrive successfully later in their study, particularly during the {\color{black} undergraduate} period. Some initial reactions from both the lecturers and the students show that both parties respond positively toward the implementation of the new curriculum.

\section*{\large Acknowledgement}
{\small We would like to thank our colleagues for composing the module descriptions and many valuable suggestions in improving the mathematics curriculum for the Foundation programme in Engineering: Lim Poay Hoon (The University of Nottingham, University Park Campus, UK), Thong Lee Fah, Koh Yew Meng (Iowa State University, USA), Grace Yap, Mohd.\ Rafi bin Segi Rahmat Harikrishan Kanthen, Balrama Applanaidu and Teo Lee Peng. Many constructive suggestions from Hong Kian Sam (University of Malaysia Sarawak), Frits van Beckum (Universiteit Twente, The Netherlands) {\color{black} and anonymous referees} are also appreciated.}

\end{document}